\documentstyle[12pt]{article}
\begin{document}

\def\dsp{\displaystyle}
\def\Rr{{bf R}}
\def\Zz{bf Z}
\def\Nn{bf N}
\def\get{\hbox{{\goth g}$^*$}}
\def\g{\gamma}
\def\om{\omega}
\def\r{\rho}
\def\a{\alpha}
\def\s{\sigma}
\def\vfi{\varphi}
\def\l{\lambda}
\def\implique{\Rightarrow}
\def\o{{\circ}}
\def\Diff{\hbox{\rm Diff}}
\def\S1{\hbox{\rm S$^1$}}
\def\Hom{\hbox{\rm Hom}}
\def\Vect{\hbox{\rm Vect}}
\def\const{\hbox{\rm const}}
\def\ad{\hbox{\hbox{\rm ad}}}
\def\semid{\hbox{\bb o}}
\def\blanc{\hbox{\ \ }}

\def\pds#1,#2{\langle #1\mid #2\rangle} %%% PRODUIT SCALAIRE
\def\f#1,#2,#3{#1\colon#2\to#3} %%%  F:A->B

\def\hfl#1{{\buildrel{#1}\over{\hbox to
12mm{\rightarrowfill}}}}

\title{Theorem on six vertices of a plane curve via the Sturm theory}

\author{L. Guieu\\
{\small Universit\'e d'Aix-Marseille III and CPT-CNRS}
\and
E. Mourre,
V.Yu. Ovsienko\\
{\small CNRS, Centre de Physique Th\'eorique}
\thanks{CPT-CNRS, Luminy Case 907,
F--13288 Marseille, Cedex 9
FRANCE
}
}

\date{}
\maketitle

\abstract{We discuss the theorem on the existence of
six points on a convex closed plane
curve in which the curve has a contact of order six
with the osculating conic. (This is the ``projective version'' of the
well known four vertices theorem for a curve in the Euclidean plane.)
We obtain this classical fact as a corollary of
some general Sturm-type theorems.}

\vfill\eject

\section{introduction}

The well known classical theorem states that a convex curve
on the Euclidean plane has
at least four vertices (critical points of its curvature). This theorem
has been frequently discussed in mathematical literature (see
\cite{arn1,tab}). Beautiful applications of this theorem to symplectic
geometry were discovered by V.I. Arnol'd \cite{arn1,arn2,arn3}. The
relation to the Sturm theory is given by S. Tabachnikov
\cite{tab}. His proof of the four vertices theorem is based on the fact
that a function on $S^1$ without $n$ first harmonics
of the Fourier decomposition vanishes at least
$2n$ times.

A point on a locally convex plane curve $c$ is called {\it sextactic}
if the oscullating conic has a
contact of order $\geq 6$ with $c$ in this point.
(Recall, that in a generic point the contact is of order 5 since a conic
is defined by 5 points.)
Sextactic points can be defined also as critical points of
the affine curvature or by the fact that the projective length element
of curve $c$ vanishes in these points.

Sextactic points are invariant under projective transformations.
This kind of singular points is an analogue of vertices
in projective (or affine) geometry. (Recall that in the
Euclidean case the osculating circle
has a contact of order $\geq 4$ with the curve in
any vertex.)

The following classical theorem can be considered as
the ``projective analogue'' of the four vertices theorem.

\proclaim Six vertices theorem. A closed convex curve on
${\bf R}^2$ has at least six sextactic points.\par

\proclaim Corollary. The affine curvature of a convex closed curve on
${\bf R}^2$ has at least six critical points.\par

The proof can be found in \cite{bol}.

The main result of this paper is a series of
general Sturm-type theorems (in the spirit of Tabachnikov
theorem). We estimate the number of zero points of
a function on $S^1$ orthogonal to solutions
of a disconjugate linear differential equation.
This approach contains at the same time the four vertices
theorem of Euclidean
geometry and the six vertices theorem.

\section{The Sturm theorems}

Consider a linear differential equation on $S^1$:
\begin{equation}
A\phi(x) =\phi^{(n)}(x)+
u_{n-1}(x)\phi^{(n-1)}(x)+
\cdots+
u_{0}(x)\phi(x)=
0
\label{equ}
\end{equation}
Here $u_i(x)\in C^{\infty}(S^1)$ (this means, all the potentials $u_i$
are periodic: $u_i(x+1)=u_i(x)$).

\vskip 0,3cm

{\bf Definition}. Equation (\ref{equ}) is called
{\it disconjugate on $S^1$} if:

{\bf 1}. Order $n=2k+1$:
all the solutions are periodic: $\phi(x+1)=\phi(x)$
and
have at most $2k$ zeros (with multiplicity) on $S^1$.

{\bf 2}. Order $n=2k$
all the solutions are anti-periodic:
$\phi(x+1)=-\phi(x)$
and
have at most $2k-1$ zeros (with multiplicity) on $S^1$.

In these cases,
$A$ is called a {\it disconjugate operator}.

\proclaim {2.1. Theorem 1}. Given a function $f\in C^{\infty}(S^1)$
orthogonal to all the solutions of a $2n+1$-order disconjugate equation:
$$
\int_{S^1}f(x)\phi(x)dx=0
$$
then $f$ has at least
$2n+2$ distinct zero points on $S^1$.\par

This is a generalization of the Tabachnikov theorem \cite{tab} stating
the same fact for a function $f(x)\in C^{\infty}(S^1)$ without $n$ first
harmonics. Indeed, such a function is orthogonal to the solutions of the
equation
$\partial_x(\partial_x^2+1)(\partial_x^2+4)\cdots
(\partial_x^2+n^2)\phi=0$.
\vskip 0,3cm
{\bf Corollary}. A function $f$ in the image of a
$2n+1$-order disconjugate operator $A$ ($f=Ag$ where
$g\in C^{\infty}(S^1)$ is any function) vanishes at least $2n+2$ times
on $S^1$.
\vskip 0,3cm
Indeed, $f$ is orthogonal to the solutions of the equation $A^*\phi=0$,
where $A^*$ is the operator adjoint to $A$. It is sufficient to remark
that the operator $A^*$ is disconjugate if $A$ is disconjugate.

\proclaim Theorem 2. Given a function $f\in C^{\infty}(S^1)$
orthogonal to all the products of solutions of a $n$-order disconjugate
equation:
$$
\int_{S^1}f(x)\phi_1(x)\phi_2(x)dx=0
$$
then $f$ has et least
$2n$ distinct zero points on $S^1$.\par

{\bf Remark}. There exist straightforward generalizations of Theorems 1 and 2.
It is sufficient
to consider a function orthogonal to a product of any 3, 4 etc.
solutions of a disconjugate equation.

\vskip 0,3cm

{\bf 2.2. Proof of Theorem 1}. Consider a function $f\in C^{\infty}(S^1)$
orthogonal
to all the solutions of a disconjugate equation $A\phi =0$.

First, observe that $f$ has at least one zero. Indeed, there
exists a solution $\phi$ positive almost everywhere on $S^1$ (take for
example a solution vanishing in some point with order $2n$, then the
disconjugacy condition implies that it has no more zero points). From
$\int_{S^1}f(x)\phi(x)dx=0$ one concludes that function $f$ changes
its sign at least once.

Let us prove that the number of points of $S^1$ in which
function $f$ has odd-order zeros (changes its sign) is superior to
$2n$. Suppose that $f$ has $2k$ odd-order zero points  $x_1,\ldots
,x_{2k}$ on $S^1$ and $k\leq n$. Consider a solution $\phi $ with two
properties:

a) $\phi $ has a zero of order $2(n-k)+1$ in $x_1$,

b) $\phi $ vanishes in all points
$x_1,\ldots x_{2k}$.

The existence of such a solution is evident. In fact, there exists a
$2k-1$-dimensional space of solutions vanishing
with order $2(n-k)+1$ in $x_1$. The subspace of this space
which consists of solutions vanishing in
$x_2$ has the dimension $\geq 2k-2$, etc.
 Now,
the disconjugacy condition implies that

a) Points $x_1,\ldots x_{2k}$ are simple zeros of $\phi $,

b) $\phi $ has no more zeros on
$S^1$.

Finally, (replacing if necessary $\phi $ by $-\phi $) one obtains
that functions $f(x)$ and $\phi (x)$ have the same sign sequence on
the segments $\rbrack x_1,x_2\lbrack,\rbrack x_2,x_3\lbrack,\ldots\rbrack
x_{2k},x_1\lbrack$ which implies the contradiction:
\hfill\break
$\int_{S^1}f(x)\phi
(x)dx>0$. The theorem is proven.

\vskip 0,3cm

{\bf 2.3. Proof of Theorem 2} is analogue to those of Theorem 1.
Suppose that $f$ has $2k$ odd-order zero points  $x_1,\ldots
,x_{2k}$ on $S^1$ and $k\leq n-1$. Take any number
$s$ which is even if $n$ is odd and odd if $n$ is even, such that
$k\leq s\leq n-1$. Then there exists a solution $\phi_1$
having odd order zero points in $x_1,\ldots ,x_s$ and such
that it has no more zero points on $S^1$ (see above). In the same way,
there exists a solution $\phi_2$
having odd order zero points in $x_{s+1},\ldots ,x_{2k}$ and such
that it has no more zero points on $S^1$. Their product
$\phi_1\phi_2$ has the same sign sequence as $f$ on
the segments $\rbrack x_1,x_2\lbrack,\rbrack x_2,x_3\lbrack,\ldots\rbrack
x_{2k},x_1\lbrack$ which implies the contradiction:
\hfill\break
$\int_{S^1}f(x)\phi_1\phi_2(x)dx\not =0$. The theorem is proven.

\section{Affine and projective lengths; affine and projective curvatures}

We recall some classical
definitions of affine and projective geometry of curves. It is very
interesting to compare the notion of length in the Euclidean, affine and
projective cases. If in the Euclidean case it measures in some sense the
distance between the curve and a fixed point, then the affine and the
projective lengths measure respectively: the distance between the curve
and a straight line, and the distance between the curve and a conic.

\vskip 0,3cm

{\bf 3.1. Affine length.} Consider a parametrised {\it locally
convex} curve $c(x)=(c_1(x),c_2(x))$ in ${\bf R}^2$ ( a curve without
inflection points). For any $x$, vectors $c'(x)$ and $c''(x)$ are
linearly independent. Define the element of affine length by
$$
d\sigma=\left|\begin{array}{cc}
c'_1(x)  & c'_2(x)\\
\noalign{\vskip 3mm}
c''_1(x) & c''_2(x)\end{array}\right|^\frac{1}{3} dx
$$
Then $\sigma $ is called the {\it affine parameter}.

\vskip 0,3cm

{\bf 3.2. Affine curvature.}
Vector $c'''(x)$ is a linear combination
of $c'(x)$ and $c''(x)$: $c'''(x)=a(x)c''(x)+b(x)c'(x)$. Moreover, the
affine parameter $\sigma $ is characterized by the fact that
$c'''(\sigma )$ is collinear to $c'$:
\begin{equation}
c'''(\sigma )=k(\sigma )c'(\sigma )
\label{aff}
\end{equation}
Function
$k(\sigma )$ is called the {\it affine curvature}.

\vskip 0,3cm

\proclaim {\bf 3.3.} Wilczynski-Cartan construction \cite{car},
\cite{wil}  (see also \cite {go}). \hfill\break
(i) A parametrised locally convex curve $c(x)\in {\bf RP}^2$
canonically defines a linear differential equation of the form:
\begin{equation}
\phi'''(x)=\kappa (x)\phi'(x)+v(x)\phi (x)
\label{pro}
\end{equation}
(ii) Any equation (\ref{pro}) uniquely defines a locally convex curve
$c(x)\subset {\bf RP}^2$
(modulo projective transformations of ${\bf RP}^2$).\par

{\bf Proof}. To associate a locally convex curve with an
equation (\ref{pro}), consider space $E$ of solutions of
(\ref{pro}). Let $V_x\subset E$ consists of solutions vanishing at the
moment $x$. One has a family of 2-dimensional subspaces in a
3-dimensional linear space, or in other words, a curve in ${\bf RP}^2$.
It is locally convex (which is easy to verify). In homogeneous coordinates,
$c=(\phi _1(x):\phi _2(x):\phi _3(x))$ where
$\phi _1(x),\phi _2(x),\phi _3(x)$ are any linearly independent
solutions  of (\ref{pro}). Therefore, the equation (\ref{pro}) is
uniquely defined by the corresponding curve.

\proclaim Lemma 1. The equation (\ref{pro}) corresponding to
a closed convex curve is disconjugate.\par

{\bf Proof}. Consider a closed convex curve $c\subset{\bf RP}^2$
(see fig.1). Such a curve
has at most two points of intersection with any projective line ${\bf
RP}^1\subset{\bf RP}^2$. In homogeneous coordinates
$c=(\phi _1(x):\phi _2(x):\phi _3(x))$ where $\phi _1(x),\phi _2(x),
\phi_3(x)$ are solutions of the corresponding equation (\ref{pro}) (see
Sec. 2.3). Therefore, any solution of (\ref{pro}) is periodic and has at
most 2 zeros on $S^1$.

\vskip 0,3cm

{\bf 3.4. Projective length.} Rewrite (\ref{pro}) in more symmetric
form:
\begin{equation}
\phi'''(x)={1\over 2}\big[\kappa (x)\phi'(x)+
(\kappa(x)\phi(x))'\big]+
h(x)\phi (x)
\label{nor}
\end{equation}
where $h(x)=v(x)-\kappa'(x)/2$. Remark here that the operator
$
A_0=\partial_x^3-
{1\over 2}(\kappa(x)\partial_x+\partial_x\kappa(x))
$
is antisymmetric.

\proclaim Definition \cite{car}. The 1-form on $c$
$
d\sigma =h(x)^{1\over3}dx
$
is called the projective length element.\par

{\bf Remark}. The quantity $h(x)$ transforms
as a cubic differential $h(x)(dx)^3$
by coordinate transformations. Therefore, the 1-form
$d\sigma $ is well defined (see \cite{go}).

The projective length shows how much the curve differs from a conic.
\proclaim Lemma 2. $c$ is a conic if and only if $h\equiv 0$.\par
{\bf Proof}. Consider a second order equation
$$
\psi ''(x)={\kappa(x)\over4}\psi (x)
$$
Verify that the solutions of
the equation $A_0\phi =0$ are given by quadratic polynomials in its
solutions. In particular, $\phi _1=\psi _1^2,\phi _2=\psi _1\psi _2,\phi
_3=\psi _2^2$ (where $\psi _1,\psi _2$ are linearly independent) is a
basis of solutions. Thus, $\phi _2^2=\phi _1\phi _3$ and the curve
$c=(\phi _1(x):\phi _2(x):\phi _3(x))$ is a conic.

\vskip 0,3cm

{\bf 3.5. Projective curvature}.
Let us suppose that $d\sigma \not =0$ and so $\sigma $ defines a local
parameter on $c$. Then, the function $\kappa(\sigma )/4$ is called the
{\it projective curvature} of the curve $c(x)$.

\vskip 0,3cm

{\bf 3.6. An affine curve as a projective curve}.
Consider a standard embedding ${\bf R}^2\hookrightarrow
{\bf RP}^2$ preserving the projective structure on ${\bf R}^2$ (see fig.2).
An affine locally convex curve $c\subset {\bf R}^2$ is embedded to
${\bf RP}^2$ as
a projective locally convex curve. To define its projective length and
projective curvature, represent the equation (\ref{aff}) in the
form (\ref{nor}):
$$
c'''(\sigma)=
{1\over 2}[k(\sigma)c'(\sigma)+(k(\sigma)c(\sigma))']-
{1\over 2}k'(\sigma)c(\sigma)
$$
Therefore, the projective length of $c$ can be defined by the
relation:
$$
h(\sigma)=-{1\over 2}k'(\sigma ).
$$
On the other hand, any projective curve can be considered (locally) as an
affine curve. The equation (\ref{pro}) reduces to the form (\ref{aff})
by a changing of the parameter.

\section{Sextactic points}
\proclaim Definition. A point of a locally convex curve $c\subset {\bf
RP}^2$ is called {\it sextactic} if there exists a conic in ${\bf RP}^2$
which has a contact of order $\geq 6$ with $c$ in this point.\par

{\bf 4.1. Critical points of the projective length}. The notion of a
sextactic point
 can be expressed in terms of the curvature (in
affine case) and in terms of the length element (in projective case).

\proclaim Proposition 1. A point of a locally
convex affine curve $c\subset {\bf R}^2$ is sextactic
if and only if it is a critical point of the affine curvature.\par

\proclaim Corollary. A point of a locally convex curve
$c\subset {\bf RP}^2$ is sextactic if and only if the
projective length element $d\sigma $ vanishes at this point.\par

Remark here that this statement is just an infinitesimal version of
Lemma 2.

\vskip 0,3cm

{\bf Proof of the proposition}. Consider a locally convex curve \break
$c\subset {\bf RP}^2$. Take the affine parameter on $c$,
then the coordinates of $c$ satisfies the equation
(\ref{aff}). In the neighborhood
of point $c_0=c(0)$ curve $c$ is given by the Taylor series:
$$
c(\sigma)=\sigma c'_0+
\frac{\sigma^2}{2}c''_0+
\frac{\sigma^3}{6}c'''_0+
\frac{\sigma^4}{24}c'^{\vee}_0+
\frac{\sigma^5}{120}c^{\vee}_0+...
$$
{}From (\ref{aff}) one has:
$$
\begin{array}{rcl}
c'''&=&kc'\\
c'^{\vee}&=&k'c'+kc''\\
c^{\vee}&=&(k''+k^2)c'+2k'c''
\end{array}
$$
and finally
$$
c(\sigma)=(\sigma+
k_0\frac{\sigma^3}{6}+
k'_0\frac{\sigma^4}{24}+...)c'_0+
(\frac{\sigma^2}{2}+
k_0\frac{\sigma^4}{24}+
k'_0\frac{\sigma^5}{120}+...)c''_0
$$
Fix coordinates $(x,y)$ on ${\bf RP}^2$ generated respectively
by vectors $c'(0)$ and $c''(0)$ (see fig.3).

Consider the following conic:
$$
x^2-2y+k_0y^2
$$
Satisfy the coordinates of $c(\sigma)$ to function
$
F(x,y)=
x^2-2y+k_0y^2.
$
One obtains:
$$
F(x(\sigma),y(\sigma))=k'_0\frac{\sigma^5}{20}+...
$$
Thus, the conic has a contact of order 5 with $c$ and so
this is the {\it osculating conic} to $c$. If $k'_0=0$, then
the order of contact is 6. The proposition is proven.

\vskip 0,3cm

{\bf 4.2. Proof of the six vertices theorem}.
Let us show how the six vertices theorem follows from
Theorem 1.

\proclaim Lemma 3. The parameter $h(x)$ in equation (\ref{pro})
satisfies the following condition:
$$
\int_{S^1}\phi _1(x)\phi _2(x)h(x)dx=
0$$
where $\phi _1(x),\phi _2(x)$ are any two solutions of (\ref{pro}).\par
{\bf Proof}. Let $\phi (x)$ be a solution of (\ref{pro}), then
$\phi h=A_0\phi$. Lemma 3 follows now from the fact that $A_0$ is
antisymmetric. Indeed,
$$\begin{array}{l}
\int_{S^1}\phi _1(x)\phi _2(x)h(x)dx=
\int_{S^1}\phi _1(x)A_0\phi _2(x)dx=\\
-\int_{S^1}A_0(\phi _1(x))\phi _2(x)dx=-
\int_{S^1}\phi _1(x)\phi _2(x)h(x)dx=0
\end{array}$$

The six vertices theorem follows now from Theorem 2 and Lemma 1.
In fact, the function
$h(x)$ is orthogonal to all the products of solutions of a disconjugate
equation
of order 3. Thus, it has at least 6 distinct zero points
on $S^1$ (Theorem 2). Sextactic points of a locally convex curve
$c\subset{\bf RP}^2$ coincide with zero points of $h$ (Proposition 1).
One obtains, that a closed convex curve on ${\bf RP}^2$ has at least 6
distinct sextactic points. The theorem is proven.

\vskip 0,3cm

{\bf 4.3. Geometrical properties of sextactic points}.
Let us give here two geometrical descriptions of sextactic points.

{\bf A}. Any curve $c$ in general position has almost everywhere a
contact of order 5 with its osculating conic. Nondegenerate sextactic
points can be characterized by the fact that $c$ does not cross its
osculating conic in such points (see fig.4).

\vskip 0,3cm

{\bf B. Dual curves}. Let $c_1$ and $c_2$ be locally convex curves, take
any two points  $p_1\in c_1$ and $p_2\in c_2$. Then, there exists a
projective transformation $Q\in PGL(3,{\bf R})$ such that $Qp_2=p_1$ and
the curve $Qc_2$ has a contact of order $\geq 5$ with $c_1$ in $p_1$. Let
$\bar c$
 be a {\it projectively dual curve} to the curve $c$.We show
that $c$ has a contact of order $\geq 5$ with $\bar c$
in sextactic
points.

\proclaim Lemma 4. A point $p$ of a locally convex curve $c\subset {\bf
RP}^2$ is sextactic if and only if $c$ there exists a projective
isomorphism $I:{\bf RP}^{2*}\buildrel{\cong}\over\to {\bf RP}^2$ such that $c$
has a
contact of order $\geq 6$ with $I(\bar c)$ in $p$.\par
{\bf Proof}. Let $C$ be the osculating conic of a locally convex curve $c$
in a point $p$. Then the dual conic $\bar C\in {\bf RP}^{2*}$ is the
osculating conic of $\bar c$. Take an isomorphism
$I:{\bf RP}^{2*}\to {\bf RP}^2$ which maps $\bar C$ to $C$ and the
point of contact of $\bar C$ with $\bar c$ to the
point of contact of $C$ with $c$.

\vskip 0,3cm

After completion of this paper we received the preprint \cite{arn4}
containing the proof of Theorem 1 and its applications to
the theory of space curves. We also discovered unpublished results of
A. Viro who gave another proof of the six vertices theorem using
the Sturm-Tabachnikov approach.

\vskip 0,3cm

{\bf Acknowledgments}.
It is a pleasure to acknowledge fruitful discussions with
V.I. Arnol'd and S. Tabachnikov.

\hfill\eject

%%%%%%%%%%%%%%%%%%%%%%%%%%%%%%%%%%%%%%%%%%%%%%%%%%%%%%%%%%%%%%%%%%%%%%%%%%%%%%
%%%%%%%%%%%%%%%%%%%%%%%%%%%%%%%%%%%%%%%%%%%%%%%%%%%%%%%%%%%%%%%%%%%%%%%%%%%%%%

\end{document}